\documentclass{math}
\usepackage{hyperref}
\hypersetup{
colorlinks=true,
urlcolor=blue,
citecolor=blue}
\usepackage[all]{xy,xypic}
\usepackage{amsfonts,amssymb,amsmath,amsgen,amsopn,amsbsy,theorem,graphicx,epsfig}
\usepackage{eufrak,amscd,bezier,latexsym,mathrsfs,enumerate,multirow}
\usepackage[utf8]{inputenc}\usepackage[english]{babel}
\usepackage[dvipsnames]{xcolor}
\usepackage[pagewise]{lineno}
\yil{2023}
\vol{47}
\fpage{794}
\lpage{815}
\doi{3394}

\title{Clairaut Riemannian maps}

\author[Meena and Yadav]{
	\textbf{Kiran MEENA\thanks{Correspondence: kirankapishmeena@gmail.com}, Akhilesh YADAV}\\
	Department of Mathematics, Institute of Science, Banaras Hindu University, Varanasi, India%,\\ ORCID iD: https://orcid.org/0000-0003-3990-857X\\
	\\ [1.8em]
	
	\rec{31.08.2022}
	\acc{14.02.2023}
	\finv{09.03.2023}
}
\amssayisi{2010 {\itshape AMS Mathematics Subject Classification:} 53B20, 53B35, 53C25.
	\newline
	\vspace{-15mm}
	\begin{center}
		\raisebox{-17ex}[0ex][0ex]{~~ \raisebox{.5ex}[0ex][0ex]{\footnotesize  This work is licensed under a Creative Commons Attribution 4.0 International License.}}
\end{center}} 

\newcommand{\bc}{\begin{center}}
	\newcommand{\ec}{\end{center}}

\numberwithin{equation}{section}

\newtheorem{theorem}{Theorem}[section]

\newtheorem{definition}[theorem]{Definition}
\newtheorem{example}[theorem]{Example}

\newtheorem{lemma}[theorem]{Lemma}

\newtheorem{remark}[theorem]{Remark}

\renewcommand{\phi}{\varphi}

\begin{document}
	
	\maketitle
	
	\begin{abstract}
		In this paper, first we define Clairaut Riemannian map between Riemannian manifolds by using a geodesic curve on the base space and find necessary and sufficient conditions for a Riemannian map to be Clairaut with a non-trivial example. We also obtain necessary and sufficient condition for a Clairaut Riemannian map to be harmonic. Thereafter, we study Clairaut Riemannian map from Riemannian manifold to Ricci soliton with a non-trivial example. We obtain scalar curvatures of $rangeF_\ast$ and $(rangeF_\ast)^\bot$ by using Ricci soliton. Further, we obtain necessary conditions for the leaves of $rangeF_\ast$ to be almost Ricci soliton and Einstein. We also obtain necessary condition for the vector field $\dot{\beta}$ to be conformal on $rangeF_\ast$ and necessary and sufficient condition for the vector field $\dot{\beta}$ to be Killing on $(rangeF_\ast)^\bot$, where $\beta$ is a geodesic curve on the base space of Clairaut Riemannian map. Also, we obtain necessary condition for the mean curvature vector field of $rangeF_\ast$ to be constant. Finally, we introduce Clairaut anti-invariant Riemannian map from Riemannian manifold to K\"ahler manifold, and obtain necessary and sufficient condition for an anti-invariant Riemannian map to be Clairaut with a non-trivial example. Further, we find necessary condition for $rangeF_\ast$ to be minimal and totally geodesic. We also obtain necessary and sufficient condition for Clairaut anti-invariant Riemannian maps to be harmonic.
		
	\end{abstract}

	\section{Introduction}
	The geometry of Riemannian submersions has been discussed widely in \cite{Falcitelli_2004}. In 1992, Fischer introduced Riemannian map between Riemannian manifolds as a generalization of an isometric immersion and Riemannian submersion that satisfies the well-known generalized eikonal equation $\| F_\ast \|^2 = rank F$, which is a bridge between geometric optics and physical optics \cite{Fischer_1992}. Further, the geometry of Riemannian maps was investigated in \cite{Sahin_2010a, Sahin_2010b, Sahin_2011a, Sahin_2014, Sahin_book, Sahin_2017b, Sahin_2017c, Akyol_2018, Akyol_2019}.
	
	An important Clairaut's relation states that $\tilde{r} sin\theta$ is constant, where $\theta$ is the angle between the velocity vector of a geodesic and a meridian, and $\tilde{r}$ is the distance to the axis of a surface of revolution. In 1972, Bishop defined Clairaut Riemannian submersion with connected fibers and gave a necessary and sufficient condition for a Riemannian submersion to be Clairaut Riemannian submersion \cite{Bishop_1972}. Further, Clairaut submersions were studied in \cite{Allison_1996, Lee_2015, Meena_ccs}.
	In \cite{Sahin_2017b}, \c{S}ahin introduced Clairaut Riemannian map by using a geodesic curve on the total space and obtained necessary and sufficient conditions for Riemannian map to be Clairaut Riemannian map. Further, \c{S}ahin gave an open problem to find characterizations for Clairaut Riemannian maps (see \cite{Sahin_2017c}, page 165, open problem 2). In Section \ref{sec3}, we introduce a new type of Clairaut Riemannian map by using a geodesic curve on the base space and obtain necessary and sufficient conditions for a Riemannian map to be Clairaut Riemannian map.
	
	A Riemannian manifold $(N,g_2)$ is called a Ricci soliton \cite{Hamilton_1966} if there exists a smooth vector field $Z_1$ (called potential vector field) on $N$ such that $\frac{1}{2} (L_{Z_1} g_2)(X_1,Y_1)$ $+ Ric (X_1,Y_1) + \lambda g_2(X_1,Y_1) = 0$, where $L_{Z_1}$ is the Lie derivative of the metric tensor of $g_2$ with respect to $Z_1$, $Ric$ is the Ricci tensor of $(N,g_2)$, $\lambda$ is a constant function and $X_1$, $Y_1$ are arbitrary vector fields on $N$. We shall denote a Ricci soliton by $(N,g_{2},Z_1,\lambda)$. The Ricci soliton $(N,g_2,Z_1,\lambda)$ is said to be shrinking, steady or expanding accordingly as $\lambda<0$, $\lambda=0$ or $\lambda>0$, respectively. It is obvious that a trivial Ricci soliton is an Einstein manifold \cite{Besse_1987} with $Z_1$ zero or Killing (Lie derivative of metric tensor $g_2$ with respect to $Z_1$ is vanishes). Ricci soliton can be used to solve the Poincar\'e conjecture \cite{Perelman_2002}. A Ricci soliton $(N,g_2,Z_1,\lambda)$ becomes an almost Ricci soliton \cite{Pigola_2011} if the function $\lambda$ is a variable. The Ricci soliton $(N,g_2,Z_1,\lambda)$ is said to be a gradient Ricci soliton if the potential vector field $Z_1$ is the gradient of some smooth function $f$ on $N$, which is denoted by $(N,g_2,f,\lambda)$. Moreover, a non-Killing tangent vector field $Z_1$ on a Riemannian manifold $(N,g_2)$ is called conformal \cite{Deshmukh_2014} if it satisfies $L_{Z_1} g_2 = 2fg_2$, where $f$ is called the potential function of $Z_1$. The submersions and Riemannian maps from a Ricci soliton to a Riemannian manifold were studied in \cite{Meric_2019a, Meena_cst, Yadav_rt, Yadav_crt, GGupta_2022}. In \cite{Yadav_rb}, present authors introduced Riemannian map from a Riemannian manifold to a Ricci soliton. In Section \ref{sec5}, we introduce Clairaut Riemannian map from a Riemannian manifold to a Ricci soliton.
	
	In \cite{Watson_1976}, Watson studied almost Hermitian submersions. In \cite{Sahin_2014}, \c{S}ahin introduced holomorphic Riemannian map as generalization of holomorphic submersion and holomorphic submanifold. In \cite{Sahin_2010a, Sahin_2011a, Akyol_2018, Akyol_2019} invariant, anti-invariant and semi-invariant Riemannian maps were studied from a Riemannian manifold to a K\"ahler manifold. Recently, present authors introduced Clairaut invariant Riemannian map from a Riemannian manifold to a K\"ahler manifold in \cite{Yadav_turk}.
	In Section \ref{sec4}, we introduce Clairaut anti-invariant Riemannian map from a Riemannian manifold to a K\"ahler manifold.
	
	\section{Preliminaries}\label{sec2}
	In this section, we recall the notion of Riemannian map between Riemannian manifolds and give a brief review of basic facts.
	
	Let $ F:(M^m,g_1)\rightarrow (N^n,g_2)$ be a smooth map between Riemannian manifolds such that $ 0 < rankF \leq min\{m,n\}$, where $dim (M)=m$ and $dim (N)=n$. We denote the kernel space of $F_\ast$ by $\nu_p= kerF_{\ast p}$ at $ p\in M $ and consider the orthogonal complementary space $\mathcal{H}_p=(kerF_{\ast p})^\bot$ to $kerF_{\ast p}$ in $T_pM$. Then the tangent space $T_pM$ of $M$ at $ p$ has the decomposition $T_pM= (kerF_{\ast p}) \oplus(kerF_{\ast p})^\bot = \nu_p \oplus \mathcal{H}_p$. We denote the range of $F_\ast $ by $rangeF_\ast$ at $ p$ $\in M$ and consider the orthogonal complementary space $(rangeF_{\ast p})^\bot$ to  $rangeF_{\ast p}$ in the tangent space $T_{F(p)}N$ of $N$ at $F(p) \in N$. Since $ rankF\leq min \{m,n \}$, we have $(rangeF_\ast)^\bot \neq \{0\}$. Thus the tangent space $T_{F(p)}N$ of $N$ at $ F(p) \in N$ has the decomposition $T_{F(p)}N= (rangeF_{\ast p}) \oplus(rangeF_{\ast p})^\bot$. Then $F$ is called Riemannian map at $p\in M$ if the horizontal restriction $F^h_{\ast p}:(kerF_{\ast p})^\bot \rightarrow (rangeF_{\ast p})$ is a linear isometry between the spaces $((kerF_{\ast p})^\bot,g_{1(p)}|_{(kerF_{\ast p})^\bot} )$ and $(rangeF_{\ast p},g_{2(p_1)}|_{(rangeF_{\ast p})})$, where $F(p)= p_1$. In other words, $F_\ast$ satisfies
	\begin{equation}\label{eqn2.1}
		g_2(F_\ast X, F_\ast Y) =g_1(X,Y),
	\end{equation}
	for all $X,Y$ vector field tangent to $\Gamma(kerF_{\ast p})^\bot$. It follows that isometric immersions and Riemannian submersions are particular Riemannian maps with $kerF_\ast= \{0\}$ and $(rangeF_\ast)^\bot= \{0\}$, respectively. The differential map $F_\ast$ of $F$ can be viewed as a section of bundle $Hom(TM,F^{-1}TN) \rightarrow M$, where $F^{-1}TN$ is the pullback bundle whose fibers at $p\in M$ is $(F^{-1}TN)_p = T_{F(p)}N$, $p\in M$. The bundle $Hom(TM,F^{-1}TN)$ has a connection $\nabla$ induced from the Levi-Civita connection $\nabla^M$ and the pullback connection $\overset{N}{\nabla^F}$. Then the second fundamental form of $F$ is given by \cite{Nore_1987}
	\begin{equation}\label{eqn2.2}
		(\nabla F_\ast) (X,Y) = \overset{N}{\nabla_{X}^F} F_\ast Y - F_\ast({\nabla}_{X}^M Y),
	\end{equation}
	for all $ X,Y \in \Gamma(TM)$, where $\overset{N}{\nabla_{X}^F} F_\ast Y \circ F= \nabla_{F_\ast X}^N F_\ast Y$. It is known that the second fundamental form is symmetric. In \cite{Sahin_2010a} \c{S}ahin proved that $(\nabla F_\ast)(X,Y)$ has no component in $rangeF_\ast$, for all $ X,Y \in \Gamma(kerF_\ast)^\bot$. More precisely, we have 
	\begin{equation}\label{eqn2.3}
		(\nabla F_\ast) (X,Y)\in \Gamma(rangeF_\ast)^\bot.
	\end{equation}
	The tension field of $F$ is defined to be the trace of the second fundamental form of $F$, i.e. $\tau(F) = trace(\nabla F_\ast) = \sum^m_{i=1}(\nabla F_\ast) (e_i, e_i)$, where $m= dim(M)$ and $\{e_1,e_2,...,e_m\}$ is the orthonormal frame on $M$. Moreover, a map $F:(M^m,g_1)\rightarrow (N^n,g_2)$ between Riemannian manifolds is harmonic if and only if the tension field of $F$ vanishes at each point $p\in M$.
	\begin{lemma}\label{lem3.2} \cite{Sahin_2010b} Let $ F:(M^m,g_1)\rightarrow (N^n,g_2)$ be a Riemannian map between Riemannian manifolds. Then the tension field of $F$ is given by $\tau (F) = -rF_\ast (H) + (m-r) H_2$, where $r= dim(kerF_\ast)$, $(m-r) = rankF$, $H$ and $H_2$ are the mean curvature vector fields of the distribution $kerF_\ast$ and $rangeF_\ast$, respectively.
	\end{lemma}
	\begin{lemma}\label{lem2.1} 
		\cite{Sahin_2011a} Let $F: (M,g_1) \rightarrow (N,g_2)$ be a Riemannian map between Riemannian manifolds. Then $F$ is umbilical Riemannian map if and only if
		\begin{equation*}
			(\nabla F_\ast)(X,Y) = g_1(X,Y) H_2,
		\end{equation*}	
		for $X,Y \in \Gamma(kerF_\ast)^\bot$ and $H_2$ is the mean curvature vector field of $rangeF_\ast$.
	\end{lemma}
	
	\noindent For any vector field $X$ on $M$ and any section $V$ of $(rangeF_\ast)^\bot,$ we have $\nabla_X^{F\bot} V$, which is the orthogonal projection of $\nabla_X^N V$ on $(rangeF_\ast)^\bot$, where $\nabla^{F \bot}$ is linear connection on $(rangeF_\ast)^\bot$ such that $\nabla^{F \bot} g_2 = 0$.
	
	\noindent Now, for a Riemannian map $F$ we define $\mathcal{S}_V$ as (\cite{Sahin_book}, p. 188)
	\begin{equation}\label{eqn2.4}
		\nabla_{F_\ast X}^N V = -\mathcal{S}_V F_\ast X + \nabla_X^{F \bot} V,
	\end{equation}
	where $\nabla^N$ is Levi-Civita connection on $N$, $\mathcal{S}_V F_\ast X$ is the tangential component (a vector field along $F$) of $\nabla_{F_\ast X}^N V.$ Thus at $p\in M$, we have $\nabla_{F_\ast X}^N V(p) \in T_{F(p)} N$, $\mathcal{S}_V F_\ast X \in F_{\ast p} (T_p M)$ and $\nabla_{X}^{F \bot} V(p) \in (F_{\ast p}(T_p M))^\bot.$ It is easy to see that $\mathcal{S}_V F_\ast X$ is bilinear in $V$, and $F_\ast X$ at $p$ depends only on $V_p$ and $F_{\ast p} X_p.$ Hence from (\ref{eqn2.2}) and (\ref{eqn2.4}), we obtain
	\begin{equation}\label{eqn2.5}
		g_2(\mathcal{S}_V F_\ast X, F_\ast Y) = g_2(V, (\nabla F_\ast)(X,Y)),
	\end{equation}
	for $X,Y \in \Gamma(kerF_\ast)^\bot$ and $V \in \Gamma(rangeF_\ast)^\bot$,
	where $\mathcal{S}_V$ is self-adjoint operator.

	\section{Clairaut Riemannian map between Riemannian manifolds}\label{sec3}
	In this section, we define Clairaut Riemannian map between Riemannian manifolds by using a geodesic curve \cite{doCarmo_book} on the base space and investigate geometry.
	
	The notion of Clairaut Riemannian map was defined by \c{S}ahin in \cite{Sahin_2017b}. According to the definition, a Riemannian map $F: (M,g_1) \rightarrow (N,g_2)$ between Riemannian manifolds is called Clairaut Riemannian map if there is a function $\tilde{r}:M \rightarrow \mathbb{R}^{+} $ such that for every geodesic $\alpha$ on $M$, the function $(\tilde{r} \circ \alpha)sin\theta$ is constant, where, for all $t$, $\theta(t)$ is the angle between $\dot{\alpha}(t)$ and the horizontal space at $\alpha(t)$.
	%\begin{theorem} $\cite{19}$ Let $F: (M,g_1) \rightarrow (N,g_2)$ be a Riemannian map with connected fibers, then $F$ is Clairaut Riemannian map with $ \tilde{r}=e^f$ if and only if each fiber is umbilical and has mean curvature vector field $H = -gradf$, where $f$ is a smooth function on $M$ and $grad f$ is the gradient of the function $f$ with respect to $g_1$.
	%\end{theorem}
	
	Thus, the notion of Clairaut Riemannian map comes from a geodesic curve on a surface of revolution. Therefore, we are going to give a definition of Clairaut Riemannian map by using geodesic curve on the base space.
	
	\begin{definition}\label{dfn3.2} 
		A Riemannian map $F: (M,g_1) \rightarrow (N,g_2)$ between Riemannian manifolds is called Clairaut Riemannian map if there is a function $\tilde{s}:N \rightarrow \mathbb{R}^{+} $ such that for every geodesic $\beta$ on $N$, the function $(\tilde{s} \circ \beta)sin\omega(t)$ is constant, where, $F_\ast X \in \Gamma(rangeF_\ast)$ for $X \in \Gamma(ker F_\ast)^\bot$ and $V \in \Gamma(rangeF_\ast)^\bot$ are components of $\dot{\beta}(t)$, and $\omega(t)$ is the angle between $\dot{\beta}(t)$ and $V$ for all $t$.
	\end{definition}
	\textbf{Note:}
	\noindent For all $U, V \in \Gamma(rangeF_\ast)^\bot$ we define
	\begin{equation*}
		\nabla^N_U V = \mathcal{R}(\nabla^N_U V) + \nabla^{F \bot}_U V,
	\end{equation*}
	where $\mathcal{R}(\nabla^N_U V)$ and $\nabla^{F \bot}_U V$ denote $rangeF_\ast$ and $(rangeF_\ast)^\bot$ part of $\nabla^N_U V$, respectively. Therefore $(rangeF_\ast)^\bot$ is totally geodesic if and only if
	\begin{equation*}
		\nabla^N_U V = \nabla^{F \bot}_U V.
	\end{equation*}
	Note that from now, throughout the paper, we are assuming $(rangeF_\ast)^\bot$ is totally geodesic.
	
	\begin {lemma}  Let $F: (M,g_1) \rightarrow (N,g_2)$ be a Riemannian map between Riemannian manifolds and $\alpha: I \rightarrow M$ be a geodesic curve on $M$. Then the curve $\beta = F \circ \alpha$ is geodesic curve on $N$ if and only if
	\begin{equation}\label{eqn3.1}
		(\nabla F_\ast)(X,X) + \nabla_X^{F \bot} V+ \nabla_V^{F \bot} V =0,
	\end{equation}
	\begin{equation}\label{eqn3.2}
		-\mathcal{S}_V F_\ast X + F_\ast(\nabla_X^M X) + \nabla_V^N F_\ast X=0,
	\end{equation}
	where $F_\ast X \in \Gamma(rangeF_\ast), V \in \Gamma(rangeF_\ast)^\bot$ are components of $\dot{\beta}(t)$ and $\nabla^N$ is Levi-Civita connection on $N$ and $\nabla^{F \bot}$ is a linear connection on $(rangeF_\ast)^\bot$.
	\end {lemma}
	\begin{proof} Let $\alpha: I \rightarrow M$ be a geodesic on $M$ with $U(t)= \nu \dot{\alpha}(t)$ and $X(t)=  \mathcal{H} \dot{\alpha}(t)$. Let $\beta = F \circ \alpha$ be a geodesic on $N$ with $F_\ast X \in \Gamma(rangeF_\ast)$ and $V \in \Gamma(rangeF_\ast)^\bot$ are components of $\dot{\beta}(t)$.
		
		\noindent Now,
		\begin{equation*}
			\nabla_{\dot{\beta}}^N \dot{\beta} =\nabla_{F_\ast X + V}^N (F_\ast X + V),
		\end{equation*}
		which implies
		\begin{equation*}
			\nabla_{\dot{\beta}}^N \dot{\beta} = \nabla_{F_\ast X}^N F_\ast X + \nabla_{F_\ast X}^N V + \nabla_{V}^N F_\ast X + \nabla_{V}^N V.
		\end{equation*}
		Using (\ref{eqn2.4}) in above equation, we get
		\begin{equation*}
			\nabla_{\dot{\beta}}^N \dot{\beta} = \overset{N}{\nabla_X^F} F_\ast X \circ F + (-\mathcal{S}_V F_\ast X + \nabla_X^{F \bot} V) + \nabla_{V}^N F_\ast X + \nabla_{V}^N V.
		\end{equation*}
		Using (\ref{eqn2.2}) in above equation, we get
		\begin{equation}\label{eqn3.3}\small
			\begin{array}{ll}
				\nabla_{\dot{\beta}}^N \dot{\beta} =(\nabla F_\ast)(X,X) + F_\ast (\nabla_X^M X) -\mathcal{S}_V F_\ast X + \nabla_X^{F \bot} V  + \nabla_{V}^N F_\ast X + \nabla_{V}^N V.
			\end{array}
		\end{equation}
		Since $(rangeF_\ast)^\bot$ is totally geodesic, (\ref{eqn3.3}) can be written as
		\begin{equation}\label{eqn3.4}\small
			\begin{array}{ll}
				\nabla_{\dot{\beta}}^N \dot{\beta} = (\nabla F_\ast)(X,X) + F_\ast (\nabla_X^M X) -\mathcal{S}_V F_\ast X + \nabla_X^{F \bot} V  + \nabla_{V}^N F_\ast X + \nabla_{V}^{F \bot} V.
			\end{array}
		\end{equation}
		Now $\beta$ is geodesic on $N$ if and only if $\nabla_{\dot{\beta}}^N \dot{\beta} =0$. Then (\ref{eqn3.4}) implies $(\nabla F_\ast)(X,X) + F_\ast (\nabla_X^M X) -\mathcal{S}_V F_\ast X + \nabla_X^{F \bot} V  + \nabla_{V}^N F_\ast X + \nabla_{V}^{F \bot} V = 0$, which completes the proof.
	\end{proof}
	\begin{theorem}\label{thm3.2} Let $F: (M,g_1) \rightarrow (N,g_2)$ be a Riemannian map between Riemannian manifolds such that $rangeF_\ast$ is connected and $\alpha$, $\beta=F \circ \alpha$ are geodesic curves on $M$ and $N$, respectively. Then $F$ is Clairaut Riemannian map with $\tilde{s}=e^g$ if and only if any one of the following conditions holds:
		\begin{enumerate}[(i)]
			\item $\mathcal{S}_V F_\ast X = -V(g) F_\ast X$, where $F_\ast X \in \Gamma(rangeF_\ast), V \in \Gamma(rangeF_\ast)^\bot$ are components of $\dot{\beta}(t)$.
			
			\item $F$ is umbilical map, and has $H_2 = -\nabla^N g$, where $g$ is a smooth function on $N$ and $H_2$ is the mean curvature vector field of $rangeF_\ast$.
		\end{enumerate}
	\end{theorem}
	\begin{proof} 	First we prove $F$ is a Clairaut Riemannian map with $\tilde{s}=e^g$ if and only if for any geodesic $\beta: I \rightarrow N$ with tangential components $F_\ast X \in \Gamma(rangeF_\ast)$ and $V \in \Gamma(rangeF_\ast)^\bot$, $t \in I$ the equation
		\begin{equation}\label{eqn3.5a}
			\begin{array}{ll}
				g_{2 \beta(t)}(F_\ast X(t), F_\ast X(t)) g_2(\dot{\beta}(t), (\nabla^N g)) + g_2\left(\mathcal{S}_V F_\ast X(t), F_\ast X(t)\right) = 0,
			\end{array}
		\end{equation}
		is satisfied. To prove this, let $\beta$ be a geodesic on $N$ with $\dot{\beta}(t) = F_\ast X(t) + V(t)$ and let $\omega(t) \in [0,\pi]$ denote the angle between $\dot{\beta}(t)$ and $V(t)$. 
		If $\dot{\beta}(t) \in \Gamma(rangeF_\ast)^\bot$, then we have $F_\ast X(t_0) = 0$ (i.e. (\ref{eqn3.5a}) is satisfied), which implies $\sin \omega (t) = 0$ at point $\beta(t_0)$. Thus for any function $\tilde{s}= e^g$ on $M$, $(\tilde{s}(\beta(t))) \sin \omega(t)$ identically vanishes. Therefore, the statement holds trivially in this case. Now, we consider the case $\sin \omega (t) \neq 0$, i.e. $\dot{\beta}(t)$ does not belongs only in $\Gamma(rangeF_\ast)^\bot$. Since $\beta$ is geodesic, its speed is constant $b = \|\dot{\beta}\|^2$ (say). Then
		\begin{equation}\label{eqn3.5}
			g_{2 \beta(t)} (V,V)= bcos^2\omega(t),
		\end{equation}
		\begin{equation}\label{eqn3.6}
			g_{2 \beta(t)} (F_\ast X,F_\ast X)= bsin^2\omega(t).
		\end{equation}
		Now differentiating (\ref{eqn3.6}) along $\beta$, we get
		\begin{equation}\label{eqn3.7}
			\frac{d}{dt} g_{2} (F_\ast X,F_\ast X)= 2bsin\omega cos\omega \frac{d \omega}{dt}.
		\end{equation}
		On the other hand,
		\begin{equation*}
			\frac{d}{dt} g_{2}(F_\ast X, F_\ast X) = 2g_2(\nabla_{\dot{\beta}}^N F_\ast X, F_\ast X).
		\end{equation*}
		By putting $\dot{\beta} = F_\ast X +V$ in above equation, we get
		\begin{equation*}
			\frac{d}{dt} g_{2}(F_\ast X, F_\ast X) = 2g_2(\nabla_{F_\ast X}^N F_\ast X + \nabla_{V}^N F_\ast X, F_\ast X),
		\end{equation*}
		which implies
		\begin{equation}\label{eqn3.8}
			\frac{d}{dt} g_{2}(F_\ast X, F_\ast X) = 2g_2(\overset{N}{\nabla_X^F} F_\ast X \circ F + \nabla_{V}^N F_\ast X, F_\ast X).
		\end{equation}
		Using (\ref{eqn2.2}) and (\ref{eqn3.2}) in (\ref{eqn3.8}), we get
		\begin{equation*}
			\frac{d}{dt} g_{2}(F_\ast X, F_\ast X) = 2g_2((\nabla F_\ast)(X,X) + F_\ast (\nabla_X^M X) +\mathcal{S}_V F_\ast X - F_\ast (\nabla_X^M X), F_\ast X).
		\end{equation*}
		Using (\ref{eqn2.3}) in above equation, we get
		\begin{equation}\label{eqn3.9}
			\frac{d}{dt} g_{2}(F_\ast X, F_\ast X) = 2g_2(\mathcal{S}_V F_\ast X, F_\ast X).
		\end{equation}
		Now from (\ref{eqn3.7}) and (\ref{eqn3.9}), we get
		\begin{equation}\label{eqn3.10}
			g_2(\mathcal{S}_V F_\ast X, F_\ast X) = bsin\omega cos\omega \frac{d \omega}{dt}.
		\end{equation}
		Moreover, $F$ is a Clairaut Riemannian map with $\tilde{s}=e^g$ if and only if $\frac{d}{dt}(e^{g \circ \beta} sin \omega) = 0$, that is,  $e^{g \circ \beta}sin \omega \frac{d(g\circ \beta)}{dt} + e^{g \circ \beta}cos \omega \frac{d \omega }{dt} =0$. By multiplying this with nonzero factor $bsin\omega$ and using (\ref{eqn3.6}), we get
		\begin{equation}\label{eqn3.11}
			g_{2}(F_\ast X, F_\ast X) \frac{d(g\circ \beta)}{dt}= -bsin\omega cos\omega \frac{d \omega}{dt}.
		\end{equation}
		Now from (\ref{eqn3.10}) and (\ref{eqn3.11}), we get
		\begin{equation*}
			g_2(\mathcal{S}_V F_\ast X, F_\ast X) = -g_{2}(F_\ast X, F_\ast X) \frac{d(g\circ \beta)}{dt},
		\end{equation*}
		which means
		\begin{equation}\label{eqn3.12}
			g_2(\mathcal{S}_V F_\ast X, F_\ast X)= -g_{2}(F_\ast X, F_\ast X) g_2(\nabla^N g, \dot{\beta}).
		\end{equation}
		Indeed assuming (\ref{eqn3.5a}) and considering any geodesic $\beta$ on $N$ with initial tangent vector which belongs in $\Gamma(rangeF_\ast)$, then by using $V(t_0) = 0$ in (\ref{eqn3.12}), we get $g$ is constant on $rangeF_\ast$ and since $rangeF_\ast$ is connected, $\nabla^N g \in \Gamma(rangeF_\ast)^\bot$. Then by (\ref{eqn3.12}), we get
		\begin{equation}\label{eqn3.13}
			g_2(\mathcal{S}_V F_\ast X, F_\ast X)= -g_{2}(F_\ast X, F_\ast X) g_2(\nabla^N g, V).
		\end{equation}
		Thus $\mathcal{S}_V F_\ast X = -V(g)F_\ast X$, where $V(g)$ is a smooth function on $N$, which implies the proof of $(i)$.
		Now, by using (\ref{eqn2.5}) in (\ref{eqn3.13}), we get
		\begin{equation}\label{eqn3.14}
			g_2(V, (\nabla F_\ast)(X,X)) = -g_{2}(F_\ast X, F_\ast X) g_2(\nabla^N g, V),
		\end{equation}
		for $F_\ast X \in \Gamma(rangeF_\ast)$ and $V \in \Gamma(rangeF_\ast)^\bot$. Now using (\ref{eqn2.2}) in (\ref{eqn3.14}), we get
		\begin{equation*}
			g_2(V, \overset{N}{\nabla_X^F} F_\ast X) = -g_2(\nabla^N g, V) g_{2}(F_\ast X, F_\ast X).
		\end{equation*}
		Thus by comparing, we get
		\begin{equation}\label{eqn3.15}
			\overset{N}{\nabla_X^F} F_\ast X = -(\nabla^N g) g_{2}(F_\ast X, F_\ast X).
		\end{equation}
		Taking trace of (\ref{eqn3.15}), we get
		\begin{equation}\label{eqn3.16}
			\sum_{j=r+1}^{m} \overset{N}{\nabla_{X_j}^F} F_\ast X_j = -(\nabla^N g)(m-r),
		\end{equation}
		where $\{ X_{r+1}, X_{r+2},...,X_m\}$ and $\{ F_\ast X_{r+1}, F_\ast X_{r+2},...,F_\ast X_m\}$ are orthonormal bases of $(kerF_\ast)^\bot$ and $rangeF_\ast$, respectively. 
		
		\noindent Moreover, the mean curvature vector field of $rangeF_\ast$ is defined by (\cite{Sahin_2010b}, \cite{Sahin_book} page 199)
		\begin{equation}\label{eqn3.17}
			H_2 = \frac{1}{m-r}\sum_{j=r+1}^{m} \overset{N}{\nabla_{X_j}^F} F_\ast X_j,
		\end{equation}
		where $\{X_j\}_{r+1 \leq j \leq m }$ is an orthonormal basis of $(kerF_\ast)^\bot$. Then from (\ref{eqn3.16}) and (\ref{eqn3.17}), we get 
		\begin{equation}\label{eqn3.18}
			H_2 = -\nabla^N g.
		\end{equation}
		Also, by (\ref{eqn3.14}), we get
		\begin{equation}\label{eqn3.19}
			(\nabla F_\ast)(X,X) = -g_{2}(F_\ast X, F_\ast X) (\nabla^N g).
		\end{equation}
		Since $F$ is Riemannian map, using (\ref{eqn2.1}) in (\ref{eqn3.19}), we get
		\begin{equation}\label{eqn3.20}
			(\nabla F_\ast)(X,X) = -g_1(X,X) (\nabla^N g).
		\end{equation}
		From (\ref{eqn3.18}) and (\ref{eqn3.20}), we get
		\begin{equation*}
			(\nabla F_\ast)(X,X)= g_1(X,X)H_2.
		\end{equation*}
		Thus by Lemma \ref{lem2.1} $F$ is umbilical map, which completes the proof.
	\end{proof}
	\begin{remark}
		In \cite{Sahin_2017b}, \c{S}ahin considered geodesic curve on the total manifold of a Riemannian map $F$, then by using Clairaut relation fibers of $F$ are totally umbilical. On the other hand, in Definition \ref{dfn3.2}, we considered geodesic curve on the base manifold of $F$, then by using Clairaut's relation $F$ becomes totally umbilical.
	\end{remark}
	
	\begin{theorem}  Let $F: (M^m,g_1) \rightarrow (N^n,g_2)$ be a Clairaut Riemannian map with $\tilde{s}=e^g$ between Riemannian manifolds such that $kerF_\ast$ is minimal. Then $F$ is harmonic if and only if $g$ is constant function on $N$.
	\end{theorem}
	\begin{proof} Since $H=0$, then by Lemma \ref{lem3.2} $F$ is harmonic if and only if $H_2=0$ if and only if $\nabla^N g =0$, which completes the proof.
	\end{proof}
	
	\begin{theorem}
		Let $F: (M^m,g_1) \rightarrow (N^n,g_2)$ be a Clairaut Riemannian map with $\tilde{s}=e^g$ between Riemannian manifolds. Then $N= N_{(rangeF_\ast)^\bot} \times_f N_{rangeF_\ast}$ is a twisted product manifold.
	\end{theorem}
	\begin{proof}
		By (\ref{eqn3.19}), (\ref{eqn3.20}) and Theorem \ref{thm3.2}, we have $\overset{N}{\nabla_{X}^F} F_\ast Y = g_1(X, Y) H_2$ for $X, Y \in \Gamma(kerF_\ast)^\bot$, which implies $rangeF_\ast$ is totally umbilical. Then proof follows by \cite{Ponge_1993}.
	\end{proof}
	\begin{example}\label{exclairaut} Let $M=\{(x_1,x_2) \in \mathbb{R}^{2}:x_1\neq 0, x_2\neq 0 \}$ be a Riemannian manifold with Riemannian metric $g_1=  e^{2x_2}dx_1^2 +dx_2^2$ on $M$. Let $N=\{(y_1,y_2) \in \mathbb{R}^{2}\}$ be a Riemannian manifold with Riemannian metric $g_2= e^{2x_2}dy_1^2 +dy_2^2$ on $N$. Consider a map $F : (M,g_1) \rightarrow (N,g_2)$ defined by
		\begin{equation*}
			F(x_1,x_2)= (x_1, 0).
		\end{equation*}
		Then, we get
		\begin{equation*}	
			kerF_\ast = span \{ U= e_2\}~\text{and}~	(kerF_\ast)^\bot = span \{ X=e_1\},
		\end{equation*}
		where $\Big\{ e_1 = e^{-x_2}\frac{\partial}{\partial x_1}, e_2 =\frac{\partial}{\partial x_2}\Big\}$ and $\Big\{ e_1' = e^{-x_2}\frac{\partial}{\partial y_1}, e_2' = \frac{\partial}{\partial y_2}\Big\}$ are bases on  $T_pM$ and $T_{F (p)}N$, respectively, for all $ p\in M$. By easy computations, we see that $F_\ast (X) =e_1'$ and $ g_1(X,X)= g_2(F_\ast X, F_\ast X)$ for $X \in \Gamma(kerF_\ast)^\bot.$ Thus $F$ is Riemannian map with $rangeF_\ast = span \{F_\ast (X) = e_1'\}$ and $(rangeF_\ast)^\bot= span \{ e_2'\}$. Now to show $F$ is Clairaut Riemannian map we will verify Theorem \ref{thm3.2}, for this we will verify (\ref{eqn3.13}). Since $V$ and $(\nabla F_\ast)(X,X) \in \Gamma(rangeF_\ast)^\bot$,e here we can write $V = a e_2'$ and $(\nabla F_\ast)(X,X)= b e_2'$ for some $a,b \in \mathbb{R}$. Then we get
		\begin{equation}\label{eqn3.21}
			g_2(V, (\nabla F_\ast)(X,X)) = g_2(a e_2', b e_2')= ab,
		\end{equation}
		and
		\begin{equation}\label{eqn3.32}
			g_2(F_\ast X, F_\ast X) = g_2(e_1', e_1') = 1.
		\end{equation}
		Since $\nabla^N g=\underset{i,j=1}{\overset{2}{\sum}} g_2^{ij} \frac{\partial g}{\partial y_i} \frac{\partial}{\partial y_j}$. Therefore for the function $g= -b y_2$
		\begin{equation}\label{eqn3.33}
			g_2(\nabla^N g, V)= -ab.
		\end{equation}
		Thus by using (\ref{eqn2.5}), (\ref{eqn3.21}), (\ref{eqn3.32}) and (\ref{eqn3.33}) we see that (\ref{eqn3.13}) holds. Thus $F$ is a Clairaut Riemannian map.
	\end{example}
	
	\section{Clairaut Riemannian map from Riemannian manifold to Ricci soliton}\label{sec5}
	In this section, we study Clairaut Riemannian map $ F:(M,g_1) \rightarrow (N,g_2)$ from a Riemannian manifold to a Ricci soliton and give some characterizations.
	\begin{lemma}
		\cite{Yadav_rb}
		Let $F: (M^m,g_1) \rightarrow (N^n,g_2)$ be a Riemannian map between Riemannian manifolds. Then the Ricci tensor on $(N,g_2)$  given by
		\begin{equation}\label{eqn5.1a}\small
			\begin{array}{ll}
				Ric(F_\ast X, F_\ast Y)&=  Ric^{rangeF_\ast}(F_\ast X, F_\ast Y) - \sum\limits_{k=1}^{n_1} \Big \{ g_2(\mathcal{S}_{\nabla_{e_k}^{F \bot} e_k} F_\ast X, F_\ast Y) \\&-  g_2(\nabla_{e_k}^N \mathcal{S}_{e_k} F_\ast X, F_\ast Y)+g_2(\mathcal{S}_{e_k}  F_\ast X, \mathcal{S}_{e_k} F_\ast Y)  +  g_2(\nabla_{e_k}^N  F_\ast X, \mathcal{S}_{e_k} F_\ast Y) \Big \},
			\end{array}
		\end{equation}
		\begin{equation}\label{eqn5.2b}\small
			\begin{array}{ll}
				Ric(V, W)&=  Ric^{(rangeF_\ast)^\bot}(V,W) - \sum\limits_{j=r+1}^{m} \Big \{ g_2(\mathcal{S}_{\nabla_{V}^{F \bot} W} F_\ast X_j, F_\ast X_j) \\&+  g_2(\mathcal{S}_{V} F_\ast X_j, \mathcal{S}_{W} F_\ast X_j)- \nabla_V^N (g_2(\mathcal{S}_{W}  F_\ast X_j, F_\ast X_j))  +  2 g_2( \mathcal{S}_{W} F_\ast X_j, \nabla_{V}^N  F_\ast X_j) \Big \},
			\end{array}
		\end{equation}
		
		\begin{equation}\label{eqn5.3c}\small
			\begin{array}{ll}
				Ric(F_\ast X,V)&= \sum\limits_{j=r+1}^{m} \Big \{ g_2((\tilde{\nabla}_X \mathcal{S})_V F_\ast X_j, F_\ast X_j) - g_2((\tilde{\nabla}_{X_j} \mathcal{S})_V F_\ast X, F_\ast X_j)\Big \} - \sum\limits_{k=1}^{n_1} g_2( R^{F \bot}(F_\ast X, e_k)V, e_k),
			\end{array}
		\end{equation}
		for $X,Y \in \Gamma(kerF_\ast)^\bot$, $V,W \in \Gamma(rangeF_\ast)^\bot$ and $F_\ast X, F_\ast Y \in \Gamma(rangeF_\ast)$, where $\{F_\ast X_j\}_{r+1 \leq j \leq m}$ and $\{ e_k\}_{1\leq k \leq n_1}$ are orthonormal bases of $rangeF_\ast$ and $(rangeF_\ast)^\bot$, respectively.
	\end{lemma}
	\begin{theorem}  
		Let $F:(M^m,g_1) \rightarrow (N^n,g_2)$ be a Clairaut Riemannian map with $\tilde{s}=e^g$ between Riemannian manifolds. Then the Ricci tensor on $(N,g_2)$  given by
		\begin{equation}\label{eqn4.24}\small
			\begin{array}{ll}
				Ric(F_\ast X, F_\ast Y)&=  Ric^{rangeF_\ast}(F_\ast X, F_\ast Y)- \sum\limits_{k=1}^{n_1} (e_k(g))^2g_2(F_\ast X, F_\ast Y) \\&+ \sum\limits_{k=1}^{n_1} g_2(\nabla_{e_k}^{F \bot} e_k, \nabla^N g) g_2(F_\ast X, F_\ast Y) - \sum\limits_{k=1}^{n_1} (\nabla_{e_k}^N e_k(g)) g_2(F_\ast X, F_\ast Y),
			\end{array}
		\end{equation}
		\begin{equation}\label{eqn4.25}
			\begin{array}{ll}
				Ric(V, W)&=  Ric^{(rangeF_\ast)^\bot}(V,W) + (m-r) g_2(\nabla^N g, \nabla^{F \bot}_V W) \\&- (m-r) V(g) W(g) - (m-r) \nabla^N_V W(g),
			\end{array}
		\end{equation}
		\begin{eqnarray}\label{eqn4.26}
			Ric(F_\ast X,V)= \sum\limits_{j=r+1}^{m}  g_2((\tilde{\nabla}_X \mathcal{S})_V F_\ast X_j, F_\ast X_j) -\sum\limits_{j=r+1}^{m} g_2((\tilde{\nabla}_{X_j} \mathcal{S})_V F_\ast X_j, F_\ast X_j)  - \sum\limits_{k=1}^{n_1} g_2( R^{F \bot}(F_\ast X, e_k)V, e_k),
		\end{eqnarray}
		for $X,Y \in \Gamma(kerF_\ast)^\bot$, $V,W \in \Gamma(rangeF_\ast)^\bot$ and $F_\ast X, F_\ast Y \in \Gamma(rangeF_\ast)$, where $\{F_\ast X_j\}_{r+1 \leq j \leq m}$ and $\{ e_k\}_{1\leq k \leq n_1}$ are orthonormal bases of $rangeF_\ast$ and $(rangeF_\ast)^\bot$, respectively.	
	\end{theorem}
	\begin{proof} Using Theorem \ref{thm3.2} and (\ref{eqn3.13}) in (\ref{eqn5.1a}), we get
		\begin{equation*}
			\begin{array}{ll}
				Ric(F_\ast X, F_\ast Y)&=  Ric^{rangeF_\ast}(F_\ast X, F_\ast Y)- \sum\limits_{k=1}^{n_1} (e_k(g))^2g_2(F_\ast X, F_\ast Y) + \sum\limits_{k=1}^{n_1} g_2(\nabla_{e_k}^{F \bot} e_k, \nabla^N g) g_2(F_\ast X, F_\ast Y) \\&- \sum\limits_{k=1}^{n_1} g_2(\nabla_{e_k}^N (e_k(g) F_\ast X), F_\ast Y) + \sum\limits_{k=1}^{n_1} g_2(\nabla^N_{e_k} F_\ast X, e_k(g)F_\ast Y),
			\end{array}
		\end{equation*}
		which implies (\ref{eqn4.24}). Also using Theorem \ref{thm3.2} and (\ref{eqn3.13}) in (\ref{eqn5.2b}), we get
		\begin{equation*}\small
			\begin{array}{ll}
				Ric(V, W)&=  Ric^{(rangeF_\ast)^\bot}(V,W) + \sum\limits_{j=r+1}^{m} g_2(\nabla_{V}^{F \bot} W, \nabla^N g) g_2(F_\ast X_j, F_\ast X_j) -\sum\limits_{j=r+1}^{m}  g_2(V(g) F_\ast X_j, W(g) F_\ast X_j) \\&- \sum\limits_{j=r+1}^{m} \nabla_V^N (g_2(W(g) F_\ast X_j, F_\ast X_j))+2\sum\limits_{j=r+1}^{m} g_2(W(g) F_\ast X_j, \nabla_{V}^N  F_\ast X_j),
			\end{array}
		\end{equation*}
		which implies (\ref{eqn4.25}). Also the proof of (\ref{eqn5.3c}) and (\ref{eqn4.26}) is same.
	\end{proof}
	\begin{theorem} Let $(N,g_2,H_2,\lambda)$ be a Ricci soliton with the potential vector field $H_2 \in \Gamma(rangeF_\ast)^\bot$ and $F:(M,g_1) \rightarrow (N,g_2)$ be a Clairaut Riemannian map with $\tilde{s}=e^g$ between Riemannian manifolds. Then
		\begin{equation*}
			s^{rangeF_\ast} =-\lambda (m-r) + (m-r) \Delta g - (m-r)(m-r-2)\|\nabla^N g\|^2,
		\end{equation*}
		where $s^{rangeF_\ast}$ is the scalar curvature of $rangeF_\ast$ and $(m-r) = \dim(range F_\ast)$.
	\end{theorem}
	\begin{proof} Since $(N,g_2,H_2,\lambda)$ admit Ricci soliton with the potential vector field $H_2 \in \Gamma(rangeF_\ast)^\bot$ then, we have
		\begin{equation*}
			\frac{1}{2}(L_{H_2}g_2)(F_\ast X, F_\ast Y) + Ric(F_\ast X, F_\ast Y) + \lambda g_2(F_\ast X, F_\ast Y) = 0,
		\end{equation*}
		for $F_\ast X, F_\ast Y \in \Gamma(rangeF_\ast)$, which implies
		\begin{equation*}
			\frac{1}{2}\{g_2(\nabla_{F_\ast X}^N H_2, F_\ast Y) + g_2(\nabla_{F_\ast Y}^N H_2, F_\ast X)\}+ Ric(F_\ast X, F_\ast Y) + \lambda g_2(F_\ast X, F_\ast Y) = 0.
		\end{equation*}
		Using (\ref{eqn2.4}) in above equation, we get
		\begin{equation*}\small
			\frac{1}{2}\{g_2(-\mathcal{S}_{H_2} F_\ast X, F_\ast Y) + g_2(-\mathcal{S}_{H_2} F_\ast Y, F_\ast X)\}+ Ric(F_\ast X, F_\ast Y) + \lambda g_2(F_\ast X, F_\ast Y) = 0.
		\end{equation*}
		Since $\mathcal{S}_{H_2}$ is self-adjoint, above equation can be written as
		\begin{equation}\label{eqn5.1}
			-g_2(\mathcal{S}_{H_2} F_\ast X, F_\ast Y)+ Ric(F_\ast X, F_\ast Y) + \lambda g_2(F_\ast X, F_\ast Y) = 0.
		\end{equation}
		Using (\ref{eqn3.13}), (\ref{eqn3.18}) and (\ref{eqn4.24}) in (\ref{eqn5.1}), we get
		\begin{equation*}
			\begin{array}{ll}
				-g_2(\nabla^N g, \nabla^N g) g_2(F_\ast X, F_\ast Y) + Ric^{rangeF_\ast}(F_\ast X, F_\ast Y)- \sum\limits_{k=1}^{n_1} (e_k(g))^2g_2(F_\ast X, F_\ast Y) \\+ \sum\limits_{k=1}^{n_1} g_2(\nabla_{e_k}^{F \bot} e_k, \nabla^N g) g_2(F_\ast X, F_\ast Y) - \sum\limits_{k=1}^{n_1} \nabla_{e_k}^N e_k(g) g_2(F_\ast X, F_\ast Y) + \lambda g_2(F_\ast X, F_\ast Y) = 0,
			\end{array}
		\end{equation*}
		where $\{ e_k\}_{1\leq k \leq n_1}$ is an orthonormal basis of $(rangeF_\ast)^\bot$. This implies
		\begin{equation}\label{eq3.5.2}
			\begin{array}{ll}
				-2\|\nabla^N g\|^2 g_2(F_\ast X, F_\ast Y) + Ric^{rangeF_\ast}(F_\ast X, F_\ast Y)\\- \sum\limits_{k=1}^{n_1} g_2(e_k, \nabla_{e_k}^{N}\nabla^N g) g_2(F_\ast X, F_\ast Y) + \lambda g_2(F_\ast X, F_\ast Y) = 0.
			\end{array}
		\end{equation}	
		Taking trace of (\ref{eq3.5.2}) for $rangeF_\ast$, we get
		\begin{equation*}
			\begin{array}{ll}
				s^{rangeF_\ast}- 2 (m-r) \|\nabla^N g\|^2 - (m-r)\sum\limits_{k=1}^{n_1} g_2(\nabla_{e_k}^{N} \nabla^N g, e_k)+ \lambda (m-r)= 0.
			\end{array}
		\end{equation*}
		Using definition of Hessian form of $g$ (i.e. $H^g (X_1,Y_1) = g_2(\nabla_{X_1}^N \nabla^N g, Y_1)$ for all $X_1, Y_1 \in \Gamma(TN)$) from \cite{Falcitelli_2004} in above equation, we get
		\begin{equation}\label{eq3.5.3}
			\begin{array}{ll}\small
				s^{rangeF_\ast}+ (m-r) \{-2\|\nabla^N g\|^2 - \sum\limits_{k=1}^{n_1}H^g(e_k, e_k)+ \lambda\}= 0.
			\end{array}
		\end{equation}
		Since we know that
		\begin{equation}\label{eq3.5.4a}
			\Delta g = \sum\limits_{j = r+1}^{m} H^g(F_\ast X_j,F_\ast X_j) + \sum\limits_{k=1}^{n_1}H^g(e_k, e_k),
		\end{equation}
		where $\{F_\ast X_j\}_{r+1 \leq j \leq m}$ and $\{ e_k\}_{1\leq k \leq n_1}$ are orthonormal bases of $rangeF_\ast$ and $(rangeF_\ast)^\bot$, respectively.
		Then by using definition of Hessian form of $g$ in (\ref{eq3.5.4a}), we get
		\begin{equation}\label{eq3.5.4b}
			\Delta g = \sum\limits_{j = r+1}^{m} g_2(\nabla_{F_\ast X_j}^{N} \nabla^N g,F_\ast X_j) + \sum\limits_{k=1}^{n_1}H^g(e_k, e_k).
		\end{equation}
		Using (\ref{eqn2.4}) in (\ref{eq3.5.4b}), we get
		\begin{equation*}
			\Delta g = -\sum\limits_{j = r+1}^{m} g_2(\mathcal{S}_{\nabla^N g} F_\ast X_j,F_\ast X_j) + \sum\limits_{k=1}^{n_1}H^g(e_k, e_k).
		\end{equation*}	
		Using Theorem \ref{thm3.2} in above equation, we get
		\begin{equation}\label{eq3.5.4c}
			\Delta g - (m-r) \|\nabla^N g\|^2 = \sum\limits_{k=1}^{n_1}H^g(e_k, e_k).
		\end{equation}
		Thus (\ref{eq3.5.3}) and (\ref{eq3.5.4c}) implies the proof.
	\end{proof}
	\begin{theorem} Let $(N,g_2,H_2,\lambda)$ be a Ricci soliton with the potential vector field $H_2 \in \Gamma(rangeF_\ast)^\bot$ and $F:(M^m,g_1) \rightarrow (N^n,g_2)$ be a Clairaut Riemannian map with $\tilde{s}=e^g$ between Riemannian manifolds. Then
		\begin{equation*}
			\begin{array}{ll}
				s^{(rangeF_\ast)^\bot} =-\lambda n_1 + (m-r+1) \Delta g  - (m-r)^2 \|\nabla^N g\|^2,
			\end{array}
		\end{equation*}
		where $s^{(rangeF_\ast)^\bot}$ denotes the scalar curvature of $(rangeF_\ast)^\bot$ and $(m-r) = \dim(rangeF_\ast)$, $n_1 = \dim((rangeF_\ast)^\bot)$.
	\end{theorem}
	\begin{proof} Since $(N,g_2,H_2,\lambda)$ admit Ricci soliton with the potential vector field $H_2 \in \Gamma(rangeF_\ast)^\bot$ then, we have
		\begin{equation*}
			\frac{1}{2}(L_{H_2}g_2)(V, W) + Ric(V, W) + \lambda g_2(V, W) = 0,
		\end{equation*}
		for $V,W \in \Gamma(rangeF_\ast)^\bot$, which implies
		\begin{equation*}
			\frac{1}{2}\{g_2(\nabla_{V}^N H_2, W) + g_2(\nabla_{W}^N H_2, V)\}+ Ric(V, W) + \lambda g_2(V, W) = 0.
		\end{equation*}
		Putting $H_2 = -\nabla^N g$ in above equation, we get
		\begin{equation}\label{eqn5.4}
			-\frac{1}{2}\{g_2(\nabla_{V}^N \nabla^N g, W) + g_2(\nabla_{W}^N \nabla^N g, V)\}+ Ric(V, W) + \lambda g_2(V, W) = 0.
		\end{equation}
		Using definition of Hessian form of $g$ and (\ref{eqn4.25}) in (\ref{eqn5.4}), we get
		\begin{equation}\label{eqn5.5}
			\begin{array}{ll}
				-H^g(V, W)+ Ric^{(rangeF_\ast)^\bot}(V,W) + (m-r) g_2(\nabla^N g, \nabla^{F \bot}_V W) \\- (m-r) V(g) W(g) - (m-r) \nabla^N_V W(g) + \lambda g_2(V, W) = 0.
			\end{array}
		\end{equation}
		Taking trace of (\ref{eqn5.5}) for $(rangeF_\ast)^\bot$, we get
		\begin{equation*}
			\begin{array}{ll}
				-\sum\limits_{k=1}^{n_1}H^g(e_k, e_k) + s^{(rangeF_\ast)^\bot} + \sum\limits_{k=1}^{n_1} (m-r) g_2(\nabla^N g, \nabla^{F \bot}_{e_k} e_k) - (m-r)\sum\limits_{k=1}^{n_1} (e_k(g))^2 - (m-r)\sum\limits_{k=1}^{n_1} \nabla^N_{e_k} e_k(g) + \lambda n_1 = 0,
			\end{array}
		\end{equation*}
		where $\{ e_k\}_{1\leq k \leq n_1}$ is an orthonormal basis of $(rangeF_\ast)^\bot$, which implies
		\begin{equation*}
			\begin{array}{ll}\small
				s^{(rangeF_\ast)^\bot}+ \lambda n_1 - (m-r)\sum\limits_{k=1}^{n_1} (e_k(g))^2 - (m-r+1)\sum\limits_{k=1}^{n_1} H^g(e_k, e_k)  = 0.
			\end{array}
		\end{equation*}
		Using (\ref{eq3.5.4c}) and $(e_k (g))^2 = g_2(\nabla^N g, e_k)^2 = g_2(\nabla^N g, \nabla^N g)$ in above equation, we get the proof.
	\end{proof}
	\begin{remark}
		Since $rangeF_\ast$ and $(rangeF_\ast)^\bot$ are subbundles of $TN$, they define distributions on $N$. Then for $F_\ast X, F_\ast Y \in \Gamma(rangeF_\ast)$, we have
		\begin{align*}
			[F_\ast X, F_\ast Y]& = \nabla_{F_\ast X}^N F_\ast Y - \nabla_{F_\ast Y}^N F_\ast X\\& = \overset{N}{\nabla_X^F} F_\ast Y \circ F -\overset{N}{\nabla_Y^F} F_\ast X \circ F.
		\end{align*}
		Using (\ref{eqn2.2}) in above equation, we get
		\begin{align*}
			[F_\ast X, F_\ast Y]& =  F_\ast(\nabla_X Y) - F_\ast(\nabla_Y X) = F_\ast(\nabla_X Y - \nabla_Y X) \in \Gamma(rangeF_\ast).
		\end{align*}
		Thus $rangeF_\ast$ is an integrable distribution. Then for any point $F(p) \in N$ there exists maximal integral manifold or a leaf of $rangeF_\ast$ containing $F(p)$.
	\end{remark}
	\begin{theorem} Let $(N,g_2,F_\ast Z,\lambda)$ be a Ricci soliton with the potential vector field $F_\ast Z \in \Gamma(rangeF_\ast)$ and $F:(M,g_1) \rightarrow (N,g_2)$ be a Clairaut Riemannian map with $\tilde{s}=e^g$ between Riemannian manifolds. Then a leaf of $rangeF_\ast$ is an almost Ricci soliton. 
	\end{theorem}
	\begin{proof} Since $(N,g_2,F_\ast Z,\lambda)$ admit Ricci soliton with the potential vector field $F_\ast Z \in \Gamma(rangeF_\ast)$ then, we have
		\begin{equation}\label{eqn5.7}
			\frac{1}{2}(L_{F_\ast Z}g_2)(F_\ast X, F_\ast Y) + Ric(F_\ast X, F_\ast Y) + \lambda g_2(F_\ast X, F_\ast Y) = 0,
		\end{equation}
		for $F_\ast X, F_\ast Y,F_\ast Z \in \Gamma(rangeF_\ast)$. Using (\ref{eqn4.24}) in (\ref{eqn5.7}), we get
		\begin{equation*}
			\begin{array}{ll}
				\frac{1}{2}(L_{F_\ast Z}g_2)(F_\ast X, F_\ast Y) + Ric^{rangeF_\ast}(F_\ast X, F_\ast Y)- \sum\limits_{k=1}^{n_1} (e_k(g))^2g_2(F_\ast X, F_\ast Y) \\+ \sum\limits_{k=1}^{n_1} g_2(\nabla_{e_k}^{F \bot} e_k, \nabla^N g) g_2(F_\ast X, F_\ast Y) - \sum\limits_{k=1}^{n_1} \nabla_{e_k}^N e_k(g) g_2(F_\ast X, F_\ast Y) + \lambda g_2(F_\ast X, F_\ast Y) = 0,
			\end{array}
		\end{equation*}
		where $\{ e_k\}_{1\leq k \leq n_1}$ is an orthonormal basis of $(rangeF_\ast)^\bot$, which implies
		\begin{equation*}
			\begin{array}{ll}
				\frac{1}{2}(L_{F_\ast Z}g_2)(F_\ast X, F_\ast Y) + Ric^{rangeF_\ast}(F_\ast X, F_\ast Y)+ \tilde{\lambda} g_2(F_\ast X, F_\ast Y) = 0,
			\end{array}
		\end{equation*}
		where $\tilde{\lambda}= - \sum\limits_{k=1}^{n_1} (e_k(g))^2+ \sum\limits_{k=1}^{n_1} g_2(\nabla_{e_k}^{F \bot} e_k, \nabla^N g) - \sum\limits_{k=1}^{n_1} e_k(e_k(g)) + \lambda$ is a smooth function on $N$. Thus a leaf of $rangeF_\ast$ is an almost Ricci soliton, which completes the proof.
	\end{proof}
	\begin{theorem} Let $(N,g_2,V,\lambda)$ be a Ricci soliton with the potential vector field $V \in \Gamma(rangeF_\ast)^\bot$ and $F:(M,g_1) \rightarrow (N,g_2)$ be a Clairaut Riemannian map with $\tilde{s}=e^g$ between Riemannian manifolds. Then a leaf of $rangeF_\ast$ is an Einstein. 
	\end{theorem}
	\begin{proof} Since $(N,g_2,F_\ast Z,\lambda)$ admit Ricci soliton with the potential vector field $F_\ast Z \in \Gamma(rangeF_\ast)$ then, we have
		\begin{equation*}
			\frac{1}{2}(L_{V}g_2)(F_\ast X, F_\ast Y) + Ric(F_\ast X, F_\ast Y) + \lambda g_2(F_\ast X, F_\ast Y) = 0,
		\end{equation*}
		for $F_\ast X, F_\ast Y \in \Gamma(rangeF_\ast)$, which implies
		\begin{equation*}
			\frac{1}{2}\{g_2(\nabla_{F_\ast X}^N V, F_\ast Y) + g_2(\nabla_{F_\ast Y}^N V, F_\ast X)\}+ Ric(F_\ast X, F_\ast Y) + \lambda g_2(F_\ast X, F_\ast Y) = 0.
		\end{equation*}
		Using (\ref{eqn2.4}) in above equation, we get
		\begin{equation*}
			\frac{1}{2}\{g_2(-\mathcal{S}_{V} F_\ast X, F_\ast Y) + g_2(-\mathcal{S}_{V} F_\ast Y, F_\ast X)\}+ Ric(F_\ast X, F_\ast Y) + \lambda g_2(F_\ast X, F_\ast Y) = 0.
		\end{equation*}
		Since $\mathcal{S}_{V}$ is self-adjoint, above equation can be written as
		\begin{equation}\label{eqn5.8}
			-g_2(\mathcal{S}_{V} F_\ast X, F_\ast Y)+ Ric(F_\ast X, F_\ast Y) + \lambda g_2(F_\ast X, F_\ast Y) = 0.
		\end{equation}
		Since $F$ is Clairaut Riemannian map, using $\mathcal{S}_V F_\ast X = -V(g) F_\ast X$ and (\ref{eqn4.24}) in (\ref{eqn5.8}), we get
		\begin{equation*}
			\begin{array}{ll}
				V(g) g_2(F_\ast X, F_\ast Y) + Ric^{rangeF_\ast}(F_\ast X, F_\ast Y)- \sum\limits_{k=1}^{n_1} (e_k(g))^2g_2(F_\ast X, F_\ast Y) \\+ \sum\limits_{k=1}^{n_1} g_2(\nabla_{e_k}^{F \bot} e_k, \nabla^N g) g_2(F_\ast X, F_\ast Y) - \sum\limits_{k=1}^{n_1} \nabla_{e_k}^N e_k(g) g_2(F_\ast X, F_\ast Y) + \lambda g_2(F_\ast X, F_\ast Y) = 0,
			\end{array}
		\end{equation*}
		where $\{ e_k\}_{1\leq k \leq n_1}$ is an orthonormal basis of $(rangeF_\ast)^\bot$, which implies
		\begin{equation*}
			Ric^{rangeF_\ast}(F_\ast X, F_\ast Y)= \lambda^{'} g_2(F_\ast X, F_\ast Y),
		\end{equation*}
		where $\lambda^{'}= \sum\limits_{k=1}^{n_1} (e_k(g))^2- \sum\limits_{k=1}^{n_1} g_2(\nabla_{e_k}^{F \bot} e_k, \nabla^N g) + \sum\limits_{k=1}^{n_1} e_k(e_k(g)) - \lambda - V(g)$ is a smooth function on $N$. Thus a leaf of $rangeF_\ast$ is an Einstein, which completes the proof.
	\end{proof}
	\begin{theorem} Let $\beta$ be a geodesic curve on $N$ and $(N,g_2,\dot{\beta},\lambda)$ be a Ricci soliton with the potential vector field $\dot{\beta}\in \Gamma(TN)$. Let $F:(M,g_1) \rightarrow (N,g_2)$ be a Clairaut Riemannian map with $\tilde{s}=e^g$ from a Riemannian manifold $M$ to an Einstein manifold $N$. Then the following statements are true:
		\begin{enumerate}[(i)]
			\item $\dot{\beta}$ is a conformal vector field on $rangeF_\ast$.
			\item  $\dot{\beta}$ is Killing vector field on $(rangeF_\ast)^\bot$ if and only if $V(g)W(g) = - H^g(V,W)$ for all $V,W \in \Gamma(rangeF_\ast)^\bot$.
		\end{enumerate} 
	\end{theorem}
	\begin{proof}  Since $(N,g_2, \dot{\beta}, \lambda)$ is a Ricci soliton then, we have
		\begin{equation}\label{eqn5.9}
			\frac{1}{2}(L_{\dot{\beta}}g_2)(F_\ast X, F_\ast Y) + Ric(F_\ast X, F_\ast Y) + \lambda g_2(F_\ast X, F_\ast Y) = 0,
		\end{equation}
		for $F_\ast X, F_\ast Y \in \Gamma(rangeF_\ast)$. Using (\ref{eqn4.24}) in (\ref{eqn5.9}), we get
		\begin{equation}\label{eqn5.10}
			\begin{array}{ll}
				\frac{1}{2}(L_{\dot{\beta}}g_2)(F_\ast X, F_\ast Y) + Ric^{rangeF_\ast}(F_\ast X, F_\ast Y)- \sum\limits_{k=1}^{n_1} (e_k(g))^2g_2(F_\ast X, F_\ast Y) \\+ \sum\limits_{k=1}^{n_1} g_2(\nabla_{e_k}^{F \bot} e_k, \nabla^N g) g_2(F_\ast X, F_\ast Y) - \sum\limits_{k=1}^{n_1} \nabla_{e_k}^N e_k(g) g_2(F_\ast X, F_\ast Y) + \lambda g_2(F_\ast X, F_\ast Y) = 0,
			\end{array}
		\end{equation}
		where $\{ e_k\}_{1\leq k \leq n_1}$ is an orthonormal basis of $(rangeF_\ast)^\bot$. Since $N$ is Einstein, putting $Ric^{rangeF_\ast}(F_\ast X, F_\ast Y)= -\lambda g_2(F_\ast X, F_\ast Y)$ in (\ref{eqn5.10}), we get
		\begin{equation*}
			\begin{array}{ll}
				\frac{1}{2}(L_{\dot{\beta}}g_2)(F_\ast X, F_\ast Y) + \mu g_2(F_\ast X, F_\ast Y) = 0,
			\end{array}
		\end{equation*}
		where $\mu= - \sum\limits_{k=1}^{n_1} (e_k(g))^2+ \sum\limits_{k=1}^{n_1} g_2(\nabla_{e_k}^{F \bot} e_k, \nabla^N g) - \sum\limits_{k=1}^{n_1} e_k(e_k(g))$ is a smooth function on $N$. Thus $\dot{\beta}$ is a conformal vector field on $rangeF_\ast$. On the other hand, since $(N,g_2, \dot{\beta}, \lambda)$ is a Ricci soliton then, we have
		\begin{equation}\label{eqn5.11}
			\frac{1}{2} (L_{\dot{\beta}} g_2)(V,W) + Ric (V,W) + \lambda g_2(V,W) = 0,
		\end{equation} 
		for any $V,W \in \Gamma(rangeF_\ast)^\bot$. Using (\ref{eqn4.25}) in (\ref{eqn5.11}), we get
		\begin{equation}\label{eqn5.12}
			\begin{array}{ll}
				\frac{1}{2} (L_{\dot{\beta}} g_2)(V,W) +  Ric^{(rangeF_\ast)^\bot}(V,W) + (m-r) g_2(\nabla^N g, \nabla^{F \bot}_V W) \\- (m-r) V(g) W(g) - (m-r) \nabla^N_V W(g) + \lambda g_2(V,W) = 0.
			\end{array}
		\end{equation}
		Since $N$ is Einstein, putting $Ric^{(rangeF_\ast)^\bot}(V,W)= -\lambda g_2(V, W)$ in (\ref{eqn5.12}), we get
		\begin{equation*}
			\begin{array}{ll}
				\frac{1}{2} (L_{\dot{\beta}} g_2)(V,W) +\{ g_2(\nabla^N g, \nabla_V^{F \bot} W) - V(g) W(g) - \nabla_V^N W(g) \}(m-r)= 0.
			\end{array}
		\end{equation*}
		Then by using $\nabla_V^N W(g) = \nabla_V^N(g_2(W, \nabla^N g))= g_2(\nabla_V^N W, \nabla^N g) + H^g(V, W)= g_2(\nabla_V^{F \bot} W, \nabla^N g) + H^g(V, W)$ in above equation, we get
		$\frac{1}{2} (L_{\dot{\beta}} g_2)(V,W)= 0$ if and only if $V(g)W(g) = -H^g(V,W)$.
		This completes the proof.
	\end{proof}
	\begin{lemma}
		Let $(N,g_2,X_1,\lambda)$ be a Ricci soliton with the potential vector field $X_1 \in \Gamma(TN)$ and $F:(M^m,g_1) \rightarrow (N^n,g_2)$ be a Clairaut Riemannian map with $\tilde{s}=e^g$ between Riemannian manifolds. Then
		\begin{equation}
			s = -\lambda n,
		\end{equation}
		where $s$ denotes the scalar curvature of $N$.
	\end{lemma}
	\begin{proof} 
		The proof is similar to remark 9 of \cite{Yadav_crt}; therefore, we are omitting it.
	\end{proof}
	\begin{theorem} 
		Let $(N,g_2,-H_2,\lambda)$ be a Ricci soliton with the potential vector field $-H_2 \in \Gamma(rangeF_\ast)^\bot$ and $F:(M,g_1) \rightarrow (N,g_2)$ be a Clairaut Riemannian map with $\tilde{s}=e^g$ between Riemannian manifolds. Then following statements are true:\\$(i)$ $N$ admits a gradient Ricci soliton.\\$(ii)$ The mean curvature vector field of $rangeF_\ast$ is constant.
	\end{theorem}
	\begin{proof}  
		By similar proof as theorem 10 of \cite{Yadav_crt}, we get
		\begin{equation*}
			\Delta g =0.
		\end{equation*}
		Hence $\nabla^N (\nabla^N g) =0$, i.e. $\nabla^N H_2 =0$, which means $H_2$ is constant. This completes the proof.
	\end{proof}
	
	\begin{example}	
		The map $F: M \to N$ given in Example \ref{exclairaut} is Clairaut Riemannian map. Now, we will show that $N$ admits a Ricci soliton, i.e.
		\begin{equation}\label{eqn5.4b} 
			\frac{1}{2} (L_{Z_1} g_2)(X_1,Y_1) + Ric (X_1,Y_1) + \lambda g_2(X_1,Y_1) = 0,
		\end{equation}
		for any $X_1, Y_1, Z_1 \in \Gamma(TN)$. 
		\noindent By similar computations as example 6.1 of \cite{Yadav_rb}, we get
		\begin{equation}\label{eqn5.8b} 
			\frac{1}{2} (L_{Z_1} g_2)(X_1,Y_1) = 0,
		\end{equation}
		\begin{equation}\label{eqn5.9b}
			g_2(X_1,Y_1) = (a_1a_3 + a_2a_4),
		\end{equation}
		and
		\begin{equation}\label{eqn5.10b}
			Ric(X_1,Y_1)= a_1a_3Ric(e_1',e_1') + (a_1a_4+a_2a_3) Ric(e_1', e_2') +a_2a_4 Ric(e_2', e_2').
		\end{equation}
		By (\ref{eqn4.24}), we get
		\begin{equation*}
			\begin{array}{ll}
				Ric(e_1', e_1')=Ric^{rangeF_\ast}(e_1', e_1')-(g_2(\nabla^N g, e_2'))^2 + g_2 (\nabla_{e_2'}^{F \bot} e_2', \nabla^N g)-\nabla_{e_2'}^N(g_2(e_2', \nabla^N g)).
			\end{array}
		\end{equation*}
		Since dimension of $rangeF_\ast$ is one,  $Ric^{rangeF_\ast}(e_1', e_1')= 0$ and we have $\nabla^N g = -be_2'$ for some $b \in \mathbb{R}$. So
		\begin{equation}\label{eqn5.11b}
			Ric(e_1', e_1')= -b^2,
		\end{equation}
		By (\ref{eqn4.25}), we get
		\begin{equation*}\small
			\begin{array}{ll}
				Ric(e_2', e_2')= Ric^{(rangeF_\ast)^\bot}(e_2', e_2') + g_2(\nabla^N g, \nabla_{e_2'}^{F \bot} e_2') - e_2'(g)e_2'(g) - \nabla_{e_2'}^N (e_2'(g)).
			\end{array}
		\end{equation*}
		Since dimension of $(rangeF_\ast)^\bot$ is one, $Ric^{(rangeF_\ast)^\bot}(e_2', e_2') = 0$ and putting $\nabla^N g = -be_2'$ for some $b \in \mathbb{R}$, we get
		\begin{equation}\label{eqn5.12b}
			Ric(e_2', e_2')= -b^2.
		\end{equation}
		And by similar computation as example 6.1 of \cite{Yadav_rb}, we get
		\begin{equation}\label{eqn5.13b}
			Ric(e_1', e_2')= 0.
		\end{equation}
		Using (\ref{eqn5.11b}), (\ref{eqn5.12b}) and (\ref{eqn5.13b}) in (\ref{eqn5.10b}), we get
		\begin{equation}\label{eqn5.14b}
			Ric(X_1,Y_1)= -(a_1a_3 + a_2a_4 )b^2.
		\end{equation}
		Now, using (\ref{eqn5.8b}), (\ref{eqn5.9b}) and (\ref{eqn5.14b}) in (\ref{eqn5.4b}), we obtain that metric $g_2$ admits Ricci soliton for
		\begin{equation*}
			\begin{array}{ll}
				\lambda = b^2.
			\end{array}
		\end{equation*}	
		Since  $b \in \mathbb{R}$, for some choices of $b$ Ricci soliton $(N,g_2)$ will be expanding or steady according to $\lambda >0$ or $\lambda = 0$.
	\end{example}
	
	\section{Clairaut anti-invariant Riemannian map from Riemannian manifold to K\"ahler manifold}\label{sec4}
	In this section, we introduce Clairaut anti-invariant Riemannian map from a Riemannian manifold to a K\"ahler manifold and investigate the geometry with a non-trivial example.
	
	Let $(N,g_2)$ be an almost Hermitian manifold \cite{Yano_1984}, then $N$ admits a tensor $J$ of type (1, 1) on $N$ such that $J^2= -I$ and
	\begin{equation}\label{eqn4.1}
		g_2(JX_1,JY_1)=g_2 (X_1,Y_1),
	\end{equation}
	for all $ X_1,Y_1 \in \Gamma(TN)$.
	An almost Hermitian manifold $N$ is called K\"ahler manifold if
	\begin{equation*}
		(\nabla_{X_1}^N J)Y_1=0,
	\end{equation*} for all  $ X_1,Y_1 \in \Gamma(TN)$, where $\nabla^N$ is the Levi-Civita connection on $N$.
	
	\begin{definition} \cite{Sahin_2010a}
		Let $F: (M,g_1) \rightarrow (N,g_2)$ be a proper Riemannian map from a Riemannian manifold $M$ to an almost Hermitian manifold $N$ with almost complex structure $J$. We say that $F$ is an anti-invariant Riemannian map at $p\in M$ if $J(rangeF_{\ast p}) \subset (rangeF_{\ast p})^\bot$. If $F$ is an anti-invariant Riemannian map for every $p\in M$ then $F$ is called an anti-invariant Riemannian map.
	\end{definition}
	In this case we denote the orthogonal subbundle to $J(rangeF_\ast)$ in $(rangeF_\ast)^\bot$ by $\mu$, i.e. $ (rangeF_\ast)^\bot = J(rangeF_\ast) \oplus \mu$. For any $V\in \Gamma(rangeF_\ast)^\bot$, we have 
	\begin{equation}\label{eqn4.2}
		JV=BV + CV,
	\end{equation}
	where $BV \in \Gamma(rangeF_\ast)$ and $CV \in \mu$. Note that if $\mu=0$ then $F$ is called Lagrangian Riemannian map \cite{Tastan_2014}.

	\begin{lemma}  Let $F: (M,g_1) \rightarrow (N,g_2,J)$ be an anti-invariant Riemannian map from a Riemannian manifold $M$ to a K\"ahler manifold $N$ and $\alpha: I \rightarrow M$ be a geodesic curve on $M$. Then the curve $\beta = F \circ \alpha$ is geodesic on $N$ if and only if
		\begin{equation}\label{eqn4.4}
			-\mathcal{S}_{J F_\ast X} F_\ast X -\mathcal{S}_{CV} F_\ast X+ \nabla_V^N BV+ F_\ast(\nabla_X^M {^\ast}F_\ast BV) =0,
		\end{equation}
		\begin{equation}\label{eqn4.5}
			(\nabla F_\ast)(X,{^\ast}F_\ast BV) + \nabla_X^{F \bot} J F_\ast X+ \nabla_V^{F \bot} JF_\ast X +\nabla_X^{F \bot}CV + \nabla_V^{F \bot} CV =0,
		\end{equation}
		where $F_\ast X \in \Gamma(rangeF_\ast), V \in \Gamma(rangeF_\ast)^\bot$ are components of $\dot{\beta}(t)$ and ${^\ast F_\ast}$ is the adjoint map of $F_\ast$, and $\nabla^N$ is the Levi-Civita connection on $N$, and $\nabla^{F \bot}$ is a linear connection on $(rangeF_\ast)^\bot$.
	\end{lemma}
	\begin{proof} Let $\alpha: I \rightarrow M$ be a geodesic on $M$ and let $\beta = F \circ \alpha$ be a geodesic on $N$ with $F_\ast X \in \Gamma(rangeF_\ast)$ and $V \in \Gamma(rangeF_\ast)^\bot$ are components of $\dot{\beta}(t)$. Since $N$ is K\"ahler manifold, $\nabla_{\dot{\beta}}^N \dot{\beta} = -J \nabla_{\dot{\beta}}^N J \dot{\beta}$. Thus
		
		\begin{equation*}
			\nabla_{\dot{\beta}}^N \dot{\beta} = -J\nabla_{\dot{\beta}}^N J\dot{\beta} = -J\nabla_{F_\ast X + V}^N J(F_\ast X + V),
		\end{equation*}
		which implies
		\begin{equation}\label{eqn4.6}
			\nabla_{\dot{\beta}}^N \dot{\beta} = -J( \nabla_{F_\ast X}^N JF_\ast X + \nabla_{F_\ast X}^N JV + \nabla_{V}^N J F_\ast X + \nabla_{V}^N JV).
		\end{equation}
		Using (\ref{eqn2.4}) and (\ref{eqn4.2}) in (\ref{eqn4.6}), we get
		\begin{equation}\label{eqn4.7}
			\begin{array}{ll}
				\nabla_{\dot{\beta}}^N \dot{\beta} &= -J (-\mathcal{S}_{J F_\ast X} F_\ast X -\mathcal{S}_{CV} F_\ast X + \nabla_{V}^N BV + \nabla_{F_\ast X}^N BV \\&+ \nabla_X^{F \bot} J F_\ast X+ \nabla_V^{F \bot} J F_\ast X + \nabla_X^{F \bot} CV + \nabla_V^{F \bot} CV).
			\end{array}
		\end{equation}
		Since $\nabla^N$ is Levi-Civita connection on $N$ and $g_2(\nabla_V^N BV,U) = 0$ for any $U\in\Gamma(rangeF_\ast)^\bot$, $\nabla_V^N BV \in \Gamma(rangeF_\ast)$ and using (\ref{eqn2.2}), we get $\nabla_{F_\ast X}^N BV= \overset{N}{\nabla_X^F} BV \circ F =(\nabla F_\ast)(X, {}^\ast F_\ast BV) + F_\ast (\nabla_X^M {}^\ast F_\ast BV)$. Then by (\ref{eqn4.7}), we get
		\begin{equation*}
			\begin{array}{ll}
				\nabla_{\dot{\beta}}^N \dot{\beta} &= -J (-\mathcal{S}_{J F_\ast X} F_\ast X -\mathcal{S}_{CV} F_\ast X + \nabla_{V}^N BV + (\nabla F_\ast)(X, {}^\ast F_\ast BV) \\&+ F_\ast (\nabla_X^M {}^\ast F_\ast BV) + \nabla_X^{F \bot} J F_\ast X+ \nabla_V^{F \bot} J F_\ast X + \nabla_X^{F \bot} CV + \nabla_V^{F \bot} CV).
			\end{array}
		\end{equation*}
		Now $\beta$ is geodesic on $N\iff$ $\nabla_{\dot{\beta}}^N \dot{\beta} =0$ $\iff$$-\mathcal{S}_{J F_\ast X} F_\ast X -\mathcal{S}_{CV} F_\ast X + \nabla_{V}^N BV + (\nabla F_\ast)(X, {}^\ast F_\ast BV) + F_\ast (\nabla_X^M {}^\ast F_\ast BV) + \nabla_X^{F \bot} J F_\ast X + \nabla_V^{F \bot} J F_\ast X + \nabla_X^{F \bot} CV + \nabla_V^{F \bot} CV=0$, which completes the proof.
	\end{proof}
	
	\begin{definition} 
		An anti-invariant Riemannian map from a Riemannian manifold to a K\"ahler manifold is called Clairaut anti-invariant Riemannian map if it satisfies the condition of Clairaut Riemannian map.
	\end{definition}
	\begin{theorem}\label{thm4.7}  Let $F: (M,g_1) \rightarrow (N,g_2,J)$ be an anti-invariant Riemannian map from a Riemannian manifold $M$ to a K\"ahler manifold $N$ and $\alpha$, $\beta=F \circ \alpha$ are geodesic curves on $M$ and $N$, respectively. Then $F$ is Clairaut anti-invariant Riemannian map with $\tilde{s}=e^g$ if and only if $g_2( \mathcal{S}_{J F_\ast X} F_\ast X +\mathcal{S}_{CV} F_\ast X, BV) - g_2((\nabla F_\ast)(X, {}^\ast F_\ast BV) + \nabla_X^{F \bot} J F_\ast X+ \nabla_V^{F \bot} J F_\ast X, CV )-g_{2}(F_\ast X, F_\ast X) \frac{d(g\circ \beta)}{dt} = 0$, where $g$ is a smooth function on $N$ and $F_\ast X \in \Gamma(rangeF_\ast),$ $V \in \Gamma(rangeF_\ast)^\bot$ are components of $\dot{\beta}(t)$.
	\end{theorem}
	\begin{proof} Let $\alpha: I \rightarrow M$ be a geodesic on $M$ and let $\beta = F \circ \alpha$ be a geodesic on $N$ with $F_\ast X \in \Gamma(rangeF_\ast)$ and $V \in \Gamma(rangeF_\ast)^\bot$ are components of $\dot{\beta}(t)$ and $\omega(t)$ denote the angle in $[0,\pi]$ between $\dot{\beta}$ and $V$. Assuming $b=\|\dot{\beta}(t)\|^2$, then we get
		\begin{equation}\label{eqn4.8}
			g_{2 \beta(t)} (V,V)= bcos^2\omega(t),
		\end{equation}
		\begin{equation}\label{eqn4.9}
			g_{2 \beta(t)} (F_\ast X,F_\ast X)= bsin^2\omega(t).
		\end{equation}
		Now differentiating (\ref{eqn4.8}) along $\beta$, we get
		\begin{equation}\label{eqn4.10}
			\frac{d}{dt} g_{2} (V,V)= -2bsin\omega(t) cos\omega(t) \frac{d \omega}{dt}.
		\end{equation}
		On the other hand by (\ref{eqn4.1}), we get
		\begin{equation*}
			\frac{d}{dt} g_{2}(V, V) = \frac{d}{dt} g_{2}(JV,JV).
		\end{equation*}
		Using (\ref{eqn4.2}) in above equation, we get
		\begin{equation*}
			\frac{d}{dt} g_{2}(V, V) = \frac{d}{dt}\Big ( g_{2}(BV,BV) +g_{2}(CV,CV)\Big ),
		\end{equation*}
		which implies
		\begin{equation}\label{eqn4.11}
			\frac{d}{dt} g_{2}(V, V) = 2g_{2}(\nabla_{\dot{\beta}}^NBV,BV) +2 g_{2}(\nabla_{\dot{\beta}}^NCV,CV).
		\end{equation}
		Putting $\dot{\beta} = F_\ast X +V$ in (\ref{eqn4.11}), we get
		\begin{equation*}
			\begin{array}{ll}
				\frac{d}{dt} g_{2}(V, V) = 2g_2({\nabla_{F_\ast X}^N} BV, BV) + 2g_2(\nabla_{F_\ast X}^N CV, CV) + 2g_2(\nabla_{V}^N BV, BV) +2g_2(\nabla_{V}^{N} CV, CV).
			\end{array}
		\end{equation*}
		Since $(rangeF_\ast)^\bot$ is totally geodesic, above equation can be written as
		\begin{equation}\label{eqn4.12}
			\begin{array}{ll}
				\frac{d}{dt} g_{2}(V, V) = 2g_2(\overset{N}{\nabla_X^F} BV \circ F, BV) + 2g_2(\nabla_{F_\ast X}^N CV, CV) + 2g_2(\nabla_{V}^N BV, BV) +2g_2(\nabla_{V}^{F \bot} CV, CV).
			\end{array}
		\end{equation}
		Using (\ref{eqn2.2}), (\ref{eqn2.3}) and (\ref{eqn2.4}) in (\ref{eqn4.12}), we get
		\begin{equation}\label{eqn4.13}
			\begin{array}{ll}
				\frac{d}{dt} g_{2}(V, V) = 2g_2(F_\ast(\nabla_X^M {}^\ast F_\ast BV) +\nabla_{V}^N BV, BV) + 2g_2(\nabla_{X}^{F \bot} CV + \nabla_{V}^{F \bot} CV, CV).
			\end{array}
		\end{equation}
		Using (\ref{eqn4.4}) and (\ref{eqn4.5}) in (\ref{eqn4.13}), we get
		\begin{equation}\label{eqn4.14}\small
			\begin{array}{ll}
				\frac{d}{dt} g_{2}(V, V) = 2g_2(\mathcal{S}_{J F_\ast X} F_\ast X +\mathcal{S}_{CV} F_\ast X, BV) -2g_2\Big( (\nabla F_\ast)(X,{^\ast}F_\ast BV) + \nabla_X^{F \bot} J F_\ast X+ \nabla_V^{F \bot} JF_\ast X, CV\Big).
			\end{array}
		\end{equation}
		Now from (\ref{eqn4.10}) and (\ref{eqn4.14}), we get
		\begin{equation}\label{eqn4.15}
			\begin{array}{ll}
				g_2(\mathcal{S}_{J F_\ast X} F_\ast X +\mathcal{S}_{CV} F_\ast X, BV) -g_2\Big( (\nabla F_\ast)(X,{^\ast}F_\ast BV) + \nabla_X^{F \bot} J F_\ast X+ \nabla_V^{F \bot} JF_\ast X, CV\Big)= -bsin\omega cos\omega \frac{d \omega}{dt}.
			\end{array}
		\end{equation}
		Moreover, $F$ is a Clairaut Riemannian map with $\tilde{s}=e^g$ if and only if $\frac{d}{dt}(e^{g \circ \beta} sin \omega) = 0$, that is,  $e^{g \circ \beta}sin \omega \frac{d(g\circ \beta)}{dt} + e^{g \circ \beta}cos \omega \frac{d \omega }{dt} =0$. By multiplying this with nonzero factor $bsin\omega$ and using (\ref{eqn4.9}), we get
		\begin{equation}\label{eqn4.16}
			g_{2}(F_\ast X, F_\ast X) \frac{d(g\circ \beta)}{dt}= -bsin\omega cos\omega \frac{d \omega}{dt}.
		\end{equation}
		Thus (\ref{eqn4.15}) and (\ref{eqn4.16}) complete the proof.
	\end{proof}
	\begin{theorem}\label{thm4.8}  Let $F: (M^m,g_1) \rightarrow (N^n,g_2,J)$ be a Clairaut anti-invariant Riemannian map with $\tilde{s}=e^g$ from a Riemannian manifold $M$ to a K\"ahler manifold $N$. Then at least one of the following statement is true:\\$(i)$ $\dim(rangeF_\ast) =1$,\\$(ii)$ $g$ is constant on $J(rangeF_\ast)$, where $g$ is a smooth function on $N$.
	\end{theorem}
	\begin{proof} Since $F$ is Clairaut Riemannian map with $\tilde{s}=e^g$ then using (\ref{eqn2.2}) in (\ref{eqn3.20}), we get
		\begin{equation}\label{eqn4.17}
			\overset{N}{\nabla_X^F} F_\ast Y - F_\ast(\nabla_X^M Y) = -g_{1}(X,Y)\nabla^N g,
		\end{equation}
		for $F_\ast Y \in \Gamma(rangeF_\ast)$ and $X,Y \in \Gamma(kerF_\ast)^\bot$. Taking inner product of (\ref{eqn4.17}) with $JF_\ast Z \in \Gamma(rangeF_\ast)^\bot$, we get
		\begin{equation}\label{eqn4.18}
			g_2(\overset{N}{\nabla_X^F} F_\ast Y - F_\ast(\nabla_X^M Y), JF_\ast Z) = -g_{1}(X,Y) g_2(\nabla^N g, JF_\ast Z).
		\end{equation}
		Since $\overset{N}{\nabla^F}$ is pullback connection of the Levi-Civita connection $\nabla^N$. Therefore $\overset{N}{\nabla^F}$ is also Levi-Civita connection. Then using metric compatibility condition in (\ref{eqn4.18}), we get
		\begin{equation*}
			-g_2(\overset{N}{\nabla_X^F} JF_\ast Z, F_\ast Y)= -g_{1}(X,Y) g_2(\nabla^N g, JF_\ast Z),
		\end{equation*}
		which implies
		\begin{equation}\label{eqn4.19}
			g_2(J\overset{N}{\nabla_X^F} F_\ast Z, F_\ast Y)= g_{1}(X,Y) g_2(\nabla^N g, JF_\ast Z).
		\end{equation}
		Using (\ref{eqn4.1}) in (\ref{eqn4.19}), we get
		\begin{equation*}
			-g_2(\overset{N}{\nabla_X^F} F_\ast Z,J F_\ast Y)= g_{1}(X,Y) g_2(\nabla^N g, JF_\ast Z).
		\end{equation*}
		Using (\ref{eqn4.17}) in above equation, we get
		\begin{equation}\label{eqn4.20}
			g_1(X,Z) g_2(\nabla^N g, J F_\ast Y)= g_{1}(X,Y) g_2(\nabla^N g, JF_\ast Z).
		\end{equation}
		Now putting $X=Y$ in (\ref{eqn4.20}), we get
		\begin{equation}\label{eqn4.21}
			g_1(X,Z) g_2(\nabla^N g, J F_\ast X)= g_{1}(X,X) g_2(\nabla^N g, JF_\ast Z).
		\end{equation}
		Now interchanging $X$ and $Z$ in (\ref{eqn4.21}), we get
		\begin{equation}\label{eqn4.22}
			g_1(X,Z) g_2(\nabla^N g, J F_\ast Z)= g_{1}(Z,Z) g_2(\nabla^N g, JF_\ast X).
		\end{equation}
		From (\ref{eqn4.21}) and (\ref{eqn4.22}), we get
		\begin{equation*}
			g_2(\nabla^N g, J F_\ast X) \left( 1 -\frac{g_{1}(X,X) g_{1}(Z,Z)}{g_1(X,Z)g_1(X,Z)}\right) = 0,
		\end{equation*}
		which implies either $\dim((kerF_\ast)^\bot) =1$ or $g_2(\nabla^N g, J F_\ast X) = 0$, which means $(JF_\ast X)(g)=0$, which completes the proof.
	\end{proof}
	\begin{theorem}\label{thm4.6}
		Let $F: (M^m,g_1) \rightarrow (N^n,g_2,J)$ be a Clairaut Lagrangian Riemannian map with $\tilde{s}=e^g$ from a Riemannian manifold $M$ to a K\"ahler manifold $N$ such that $\dim(rangeF_\ast)>1$. Then following statements are true:\\
		$(i)$ $rangeF_\ast$ is minimal.\\
		$(ii)$ $rangeF_\ast$ is totally geodesic.
	\end{theorem}
	\begin{proof}
		Since $F$ is Clairaut Riemannian map then from (\ref{eqn3.20}) and Theorem \ref{thm3.2}, we have
		\begin{equation*}
			(\nabla F_\ast)(X, X) = g_1(X, X) H_2,
		\end{equation*}
		for $X\in \Gamma(kerF_\ast)^\bot$ and $H_2$ is the mean curvature vector field of $rangeF_\ast$. Now multiply above equation by $U \in \Gamma(rangeF_\ast)^\bot$, we get
		\begin{equation}\label{eqn5}
			g_2((\nabla F_\ast)(X, X), U) = g_1(X, X) g_2(H_2, U).
		\end{equation}
		Using (\ref{eqn2.2}) in (\ref{eqn5}), we get
		\begin{equation}\label{eqn6}
			g_2(\overset{N}{\nabla_{X}^F} F_\ast X, U) = g_1(X, X) g_2(H_2, U).
		\end{equation}
		Since $N$ is K\"ahler manifold, using (\ref{eqn4.1}) in (\ref{eqn6}), we get 
		\begin{equation}\label{eqn7}
			g_2(\overset{N}{\nabla_{X}^F} J F_\ast X, JU) = g_1(X, X) g_2(H_2, U).
		\end{equation}
		Since $\nabla^N$ is Levi-Civita connection on $N$, using metric compatibility condition in (\ref{eqn7}), we get
		\begin{equation}\label{eqn8}
			-g_2(J F_\ast X, \overset{N}{\nabla_{X}^F} JU) = g_1(X, X) g_2(H_2, U).
		\end{equation}
		Using (\ref{eqn6}) in (\ref{eqn8}), we get
		\begin{equation}\label{eqn9}
			-g_2(J F_\ast X, g_1(X, {}^\ast F_\ast JU)H_2) = g_1(X, X) g_2(H_2, U),
		\end{equation}
		where ${^\ast F_\ast}$ is the adjoint map of $F_\ast$. Now using $H_2 =-\nabla^N g$ in (\ref{eqn9}), we get
		\begin{equation*}
			g_1(X, {}^\ast F_\ast JU) g_2(J F_\ast X, \nabla^N g) = g_1(X, X) g_2(H_2, U),
		\end{equation*}
		which implies
		\begin{equation}\label{eqn10}
			g_1(X, {}^\ast F_\ast JU) J F_\ast X(g) = g_1(X, X) g_2(H_2, U).
		\end{equation}
		Since $\dim(rangeF_\ast)>1$ then by Theorem \ref{thm4.8}, $g$ is constant on $J(rangeF_\ast)$, which means $JF_\ast X(g) =0$. Then (\ref{eqn10}) implies
		$	g_2(H_2, U)= 0$.
		Thus
		\begin{equation}\label{eqn11}
			H_2 =0,
		\end{equation}
		which implies $(i)$.
		
		\noindent Since $H_2= trace~(\overset{N}{\nabla_{X}^F} F_\ast Y)$. Then by (\ref{eqn11}), we get $\overset{N}{\nabla_{X}^F} F_\ast Y = 0$, which implies $(ii)$.
	\end{proof}
	\begin{theorem}
		Let $F: (M^m,g_1) \rightarrow (N^n,g_2,J)$ be a Clairaut Lagrangian Riemannian map with $\tilde{s}=e^g$ from a Riemannian manifold $M$ to a K\"ahler manifold $N$ such that $\dim(rangeF_\ast)>1$. Then $F$ is harmonic if and only if mean curvature vector field of $kerF_\ast$ is constant.
	\end{theorem}
	\begin{proof}
		Let $ F:(M^m,g_1)\rightarrow (N^n,g_2)$ be a smooth map between Riemannian manifolds. Then $F$ is harmonic if and only if the tension field $\tau(F)$ of map $F$ vanishes. Then proof follows by Lemma \ref{lem3.2} and Theorem \ref{thm4.6}.
	\end{proof}
	
	\begin{theorem}
		Let $F: (M^m,g_1) \rightarrow (N^n,g_2,J)$ be a Clairaut Lagrangian Riemannian map with $\tilde{s}=e^g$ from a Riemannian manifold $M$ to a K\"ahler manifold $N$ such that $\dim(rangeF_\ast)>1$. Then $N= N_{rangeF_\ast} \times N_{(rangeF_\ast)^\bot}$ is a usual product manifold.
	\end{theorem}
	\begin{proof}
		The proof follows by \cite{Ponge_1993} and Theorem \ref{thm4.6}.
	\end{proof}
	\begin{example}  Let $M=\{(x_1,x_2) \in \mathbb{R}^{2}:x_1\neq 0, x_2\neq 0\}$ be a Riemannian manifold with Riemannian metric $g_1= e^{2x_2}dx_1^2 + e^{2x_2} dx_2^2$ on $M$. Let $N=\{(y_1,y_2) \in \mathbb{R}^{2}\}$ be a Riemannian manifold with Riemannian metric $g_2= e^{2x_2} dy_1^2 + dy_2^2$ on $N$ and the complex structure $J$ on $N$ defined as $J(y_1,y_2)=(-y_2,y_1)$. Consider a map $F : (M,g_1) \rightarrow (N,g_2,J)$ defined by
		\begin{equation*}
			F(x_1,x_2)= \Big( \frac{x_1 - x_2}{\sqrt{2}},0 \Big).
		\end{equation*}
		Then
		\begin{equation*}	
			kerF_\ast = span \Big\{ U= \frac{e_1 + e_2}{\sqrt{2}} \Big\}~\text{and}~	(kerF_\ast)^\bot = span \Big\{ X= \frac{e_1 - e_2}{\sqrt{2}} \Big\},
		\end{equation*}
		where $\Big\{ e_1=e^{-x_2}\frac{\partial}{\partial x_1},e_2 =e^{-x_2}\frac{\partial}{\partial x_2} \Big\}$ and $\Big\{ e_1' = e^{-x_2}\frac{\partial}{\partial y_1}, e_2' = \frac{\partial}{\partial y_2}\Big\}$ are bases on  $T_pM$ and $T_{F (p)}N$ respectively, for $p\in M$. By easy computations, we see that $F_\ast (X) = e_1'$ and $ g_1(X,X)= g_2(F_\ast X, F_\ast X)$ for $X \in \Gamma(kerF_\ast)^\bot.$ Thus $F$ is Riemannian map with $rangeF_\ast = span \{F_\ast (X) =e_1' \}$ and $(rangeF_\ast)^\bot= span \{ e_2'\}.$ Moreover it is easy to see that $J F_\ast X= J e_1' = -e_2'$. Thus $F$ is an anti-invariant Riemannian map.
		
		Now to show $F$ is Clairaut Riemannian map we will find a smooth function $g$ on $N$ satisfying $(\nabla F_\ast)(X,X) = -g_1(X,X) \nabla^N g$ for $X\in \Gamma(kerF_\ast)^\bot$. Since $(\nabla F_\ast)(X,X) \in \Gamma(rangeF_\ast)^\bot$ for any $X\in \Gamma(kerF_\ast)^\bot$. So here we can write $(\nabla F_\ast)(X,X) = ae_2'$, for some $a\in \mathbb{R}$. Since
		$	\nabla^N g= e^{-2x_2}\frac{\partial g}{\partial y_1} \frac{\partial}{\partial y_1} + \frac{\partial g}{\partial y_2} \frac{\partial}{\partial y_2}$.
		Hence $\nabla^N g= -a\frac {\partial}{\partial y_2}= -a e_2'$ for the function $g= -ay_2$. Then it is easy to verify that $(\nabla F_\ast)(X,X) = -g_1(X,X) \nabla^N g$, where $g_1(X,X)=1$, for vector field $X\in \Gamma(kerF_\ast)^\bot$ and we can easily see that $\nabla^N_{e_2'}e_2' = 0$. Thus by Theorem \ref{thm3.2}, $F$ is Clairaut anti-invariant Riemannian map.
	\end{example}
	\section*{Acknowledgment}
	We would like to thank all the anonymous referees for his/her valuable comments and suggestions towards the improvement of quality of the paper. The first author, Kiran Meena gratefully acknowledges the financial support provided by the Human Resource Development Group - Council of Scientific and Industrial Research (HRDG-CSIR), New Delhi, India [File No.: 09/013(0887)/2019-EMR-I].

\end{document}